\documentclass{article}
\newtheorem{thm}{Theorem}[section]

\newtheorem{prop}{Proposition}[section]

\newtheorem{lem}{Lemma}[section]

\newtheorem{ex}{Example}[section]

\newcommand{\be}{\begin{equation}}
\newcommand{\ee}{\end{equation}}
\newcommand{\beq}{\begin{eqnarray}}
\newcommand{\eeq}{\end{eqnarray}}
\newcommand{\beqn}{\begin{eqnarray*}}
\newcommand{\eeqn}{\end{eqnarray*}}
\newcommand{\e}{\varepsilon}
\newcommand{\pa}{\partial}

\newcommand{\R}{{\rm R}}

\newcommand{\tF}{\tilde{F}}

\newcommand{\pxi}{ {\pa \over \pa x^i}}

\newcommand{\pyi}{{\pa \over \pa y^i}}

\newcommand{\qed}{\hspace*{\fill}Q.E.D.}  

\title{\bf On Projectively Related\\
 Einstein Metrics}
\author{Zhongmin Shen}
\date{October 10, 1999}
\begin{document}
\maketitle

\begin{abstract}
In this paper we study pointwise projectively related Einstein metrics (having the same geodesics as point sets). 
We show that   pointwise projectively related  Einstein metrics  satisfy a  simple equation along geodesics.  
 In particular, we show that if two  pointwise projectively related  Einstein metrics  are complete with  negative Einstein constants ,  then one is a multiple of another. 
\end{abstract}

\section{Introduction}

Two regular metrics on a manifold are said to be {\it pointwise projectively related} if they have the same geodesics as point sets. 
Two regular  metric spaces  are said to be {\it projectively related} if there is a diffeomorphism between them such that the pull-back metric is pointwise projective to another one. 
Regular  metrics under our consideration are Finsler metrics neither necessarily Riemannian
nor reversible. A Riemann metric  has two features: length and angle, while a Finsler metric has only one feature: length. In projective metric geometry, the second feature of 
 Riemann metrics is generally not used.  Thus, we do not restrict 
our attention to Riemannian  metrics. 

There are many Finsler metrics on a  strongly convex subset $\Omega \subset {\rm R}^n$ which are pointwise projective to the standard Euclidean metric (such Finsler  metrics are simply called the {\it projective Finsler metrics}). 
  The  problem of characterizing and studying projective Finsler metrics is known as {\it Hilbert's fourth problem}. In dimension 2, Darboux gave a general formula for projective Finsler metrics \cite{Da}\cite{M2}.  R. Alexander's paper \cite{Ale} treats the planar case of Hilbert's fourth problem by a direct geometric argument which is entirely free from the notion of differentiability.   
In \cite{Al1}, J.C. \'{A}lvarez presents two constructions of projective metrics on ${\rm R}^n$. 
In \cite{AGS}, \'{A}lvarez-Gelfand-Smirnov presents more comprehensive  constructions of projective metrics
on an open convex subset $\Omega \subset {\rm R}^n$. See also \cite{Al2}. 
Although there are many projective metrics,  
all of them  are of scalar curvature. 
 If a projective Finsler metric is an Einstein metric, then it must be  of 
constant curvature. See Section  2 below for more details.

On an arbitrary strongly convex domain $\Omega \subset \R^n$, there is a 
 pair of Funk metrics $F_{\pm}$. The  symmetrization
$F_H : = {1\over 2} (F_{-}+F_{+})$ is called the {\it   Hilbert metric} on $\Omega$. It is known that $F_{\pm}$ are positively/negatively complete with  constant curvature $-1/4$, while $F_H$ is  complete with  constant curvature $-1$. All of them are pointwise projective to  the Euclidean metric $F_E(y) := |y|$. See \cite{Ok} for a beautiful proof on this fact.  The pair of Funk metrics on the unit ball ${\rm B}^n \subset \R^n$ are given by
\be
{F}_{\pm}(y): ={ \sqrt{  |y|^2 - \Big ( |x|^2 |y|^2 - \langle x, y \rangle^2 \Big ) }
 \pm \langle x, y \rangle 
 \over  1 -|x|^2  },\label{BaoShenMetric}
\ee
where $|\cdot |$ and $\langle \; , \; \rangle$ denotes the standard Euclidean norm and inner product.
 Note that  the Hilbert metric $F_{H} = {1\over 2} (F_{-}+F_{+})$ on ${\rm B}^n$ is just the 
 Klein metric $F_K$ given by 
\be
F_K (y): ={\sqrt{|y|^2 - \Big ( |x|^2|y|^2 - \langle x, y \rangle^2 \Big )}  \over 1- |x|^2 }.\label{ccHilbert}
\ee

\bigskip

A natural problem is to determine all  projective metrics of constant curvature
on a given open subset in $\R^n$. 
More general, 
given  a Finsler metric on a manifold $M$, we would like  to  determine 
all  Finsler  metrics which are pointwise projective to  the given one. 
There are several papers  on this problem, especially in Riemann geometry. See \cite{Mi} 
for a survey. 

In this paper, we study the following problem: given an Einstein metric,
describe all Einstein metrics which are pointwise projective to the given one. 
Below are our main theorems.

\begin{thm}\label{thm1.1} Let  $F$ and  $\tilde{F}$ be 
 Einstein metrics on a closed $n$-manifold $M$ with 
\[ {\bf Ric} =(n-1)  \lambda,  \ \ \ \ \widetilde{\bf Ric} = (n-1) \tilde{\lambda}  ,\]
where $\lambda, \tilde{\lambda} \in \{  -1, 0, 1\}$. 
Suppose that $\tilde{F}$ is pointwise projectively related to $F$. 
Then $\lambda $ and $\tilde{\lambda}$ have the same sign. More details are given below.

(i) If $\lambda =1=\tilde{\lambda}$,
then along any unit speed geodesic $c(t)$ of $F$
\be
\tilde{F}(\dot{c}(t))= {2 \over \Big ( a^2-1/a^2-b^2\Big ) \cos (2t) + 2 ab \sin (2t) +\Big ( a^2 +1/a^2 +b^2\Big ) } ,\label{lambda=1(i)}
\ee
where $a>0$ and $-\infty < b < \infty$ are constants.
Thus, for any unit speed geodesic segment $c$ of $F$ with length of $\pi$,
 it is also  a geodesic segment  of $\tilde{F}$ (as a point set) with length of $\pi$.

(ii) If $\lambda=0=\tilde{\lambda}$, then along any geodesic $c(t)$ of $F$ or $\tilde{F}$ ,
\be
 {F(\dot{c}(t)) \over \tilde{F}(\dot{c}(t))} = constant .\label{constantF}
\ee

(iii)
If $\lambda =-1=\tilde{\lambda}$, then 
\[ \tilde{F} =  {F}.\]
\end{thm}

\bigskip
Theorem \ref{thm1.1} is stated for Finslerian Einstein metrics. 
Since Riemann metrics are special Finsler metrics, this is of course  true   for Riemannian Einstein metrics.

\bigskip

 By Theorem \ref{thm1.1}(i), we know that for any Einstein metric $F$ with ${\bf Ric}=n-1$ on ${\rm S}^n$, if the geodesics of $F$ are great circles on  ${\rm S}^n$, then along any unit speed great circle $c(t)$ on ${\rm S}^n$,
${F}(\dot{c}(t))$ must be in the form (\ref{lambda=1(i)}) and the $F$-length of $c$ is still equal to $2\pi$.

According to \cite{Br1}\cite{Br2}, there are lots of non-reversible Finsler metrics 
$F$ of constant curvature $K=1$ on ${\rm S}^2$ whose geodesics are great circles
as point sets with $F$-length of $2\pi$. 
The later is also guaranteed by Theorem \ref{thm1.1}(i). In Example \ref{Bryantmetric} below, we shall extend Bryant's metrics to higher dimensional spheres.

\bigskip

Theorem \ref{thm1.1} (ii) is an almost projective rigidity result for Ricci-flat metrics. In general,  Ricci-flat  metrics are not projectively isolated. 
 For example, 
locally Minkowski metrics on a torus ${\rm T}^n$  are pointwise projective to 
the standard flat Riemann metric  on ${\rm T}^n$. In fact, these are the only flat metrics on ${\rm T}^n$.

\bigskip

Theorem \ref{thm1.1} (iii) is  a projective rigidity result for negative Einstein metrics. Any negative Einstein metric on a closed manifold is projectively isolated. 

\bigskip

\begin{thm}\label{thm1.2}
Let $F$ and $\tilde{F}$ be Ricci-flat metrics on a non-compact $n$-manifold $M$.
Suppose that $F$ and 
 $\tilde{F}$ are  pointwise projectively related. Then 
 along any unit speed geodesic $c(t)$ of $F$,
\be
\tilde{F}(\dot{c}(t)) ={1\over ( a+ bt)^2 },
\ee
where $a >0$ and $-\infty < b < \infty$ are constants.
Then  $F$ is complete if and only if  $\tilde{F}$ is complete. In this case, along any geodesic $c(t)$ of $F$ or $\tilde{F}$,
\[{F(\dot{c}(t))\over \tilde{F}(\dot{c}(t))} = constant.\]

\end{thm}

Any  Minkowski metric on ${\rm R}^n$ is pointwise projective to the standard Euclidean metric   on ${\rm R}^n$.

We actually prove that 
if  $c$ is a unit speed geodesic of $F$,  along which $ F/\tilde{F} \not=constant$,
then $c$ can not be defined on $(-\infty, \infty)$. If $c$ is defined on $[0, \infty)$, then it has finite 
$\tilde{F}$-length. If $c$ is defined on $(-\infty, 0]$, then it has finite $\tilde{F}$-length. 
This  suggests that 
there might be  a  positively
(Ricci-)flat metric
and a negatively complete (Ricci-)flat metric 
which are pointwise projective to each other (they might be also pointwise projective to the standard Euclidean metric on an open subset $\Omega \subset {\rm R}^n$). Such examples have not been found yet.

\bigskip

\begin{thm}\label{thm1.3}
Let $F$ and $\tilde{F}$ be Einstein metrics on a non-compact $n$-manifold with
\[ {\bf Ric} = -(n-1),  \ \ \ \ \widetilde{\bf Ric} = -(n-1) .\]
Suppose that
 $\tilde{F}$ is pointwise projectively related to $F$, 
 then along any geodesic $c(t)$ in $(M, F)$,
\be
\tilde{F} (\dot{c}(t)) 
= {2\over 
\Big (-1/a^2 + a^2 +b^2 \Big )\cosh (2t) 
+ 2 ab \sinh(2t) - \Big (-1/a^2 -a^2 +b^2 \Big ) },
\ee
where $ a>0$ and $-\infty < b < \infty$ are constants.

(i) If both $F$ and $\tilde{F}$ are complete. Then 
 \[ F=\tilde{F} .\]

(ii) If $F$ is complete, then for any geodesic $c$ of $F$, the $\tilde{F}$-length 
of $c$ is finite unless
\be
 \tilde{F}(\dot{c}(t)) = {1\over e^{\pm 2t} (a^2-1) + 1 } , \ \ \ \ \ (a \geq 1).\label{9a}
\ee 
\end{thm}

\bigskip

According to Theorem \ref{thm1.3}, if a Finsler metric $F$ is a complete projective metric on an open subset $\Omega \subset \R^n$ with constant curvature $-1$, then it must be the Hilbert metric 
$F_H$ on $\Omega$. This fact was proved earlier in \cite{Ber}\cite{Fk}. Theorem \ref{thm1.3}
also gives us some information on incomplete projective metrics of constant curvature $-1$.

 One can verify that for  $F:= F_{H}$ and $\tilde{F}:=F_{\pm}$, along any unit speed  geodesic $c(t)$ of $F$, 
$ \tilde{F} (\dot{c}(t))$ is in the form (\ref{9a}).
 Further discussions on examples are given in \S 7 below.

\bigskip

\section{Preliminaries}

Let $M$ be an $n$-dimensional manifold. A Finsler metric under our consideration is a function $F: TM \to [0, \infty)$ with the following properties.

(i) $ F$ is $C^{\infty}$ on $TM-\{0\}$,

(ii) For any $x\in M$, $F_x := F |_{T_xM}$ is a Minkowski functional, namely,

(iia) $F_x$ is positively homogeneous of degree one
\[ F_x(\lambda y) = \lambda F_x(y), \ \ \ \ \ \lambda >0, \ y\in T_xM\]

(iib) for any $y\in T_xM-\{0\}$, the Hessian $g_{ij}(y)$ of ${1\over 2} F^2$ at $y$ is positive definite, where
\[ g_{ij} (y) := {1\over 2} {\pa ^2 [F^2]\over \pa y^i \pa y^j } (y).\]

\bigskip
For any $y\in T_xM-\{0\}$, 
the Hessian $(g_{ij}(y))$ induces an inner product $g_y$ in $T_xM$ by
\[ g_y(u, v) : = g_{ij}(y) u^i v^j .\]
If $g$ is a Riemann metric, then 
\[ F(y):=\sqrt{g(y,y)}\]
 is a Finsler metric. 

\bigskip
Every Finsler metric induces a spray ${\bf G}$ on $M$ 
\be
{\bf G} = y^i \pxi - 2 G^i (y) \pyi,
\ee
where
\be
G^i (y) := {1\over 4} g^{il} \Big \{ {\pa^2 [F^2]\over \pa x^k \pa y^l}y^k - {\pa [F^2]\over \pa x^l } \Big \} .
\ee
${\bf G}$ is a globally defined vector field on $TM$. The projection of a flow line of ${\bf G}$ is called a {\it geodesic} 
in $M$. In local coordinates, a curve $c(t)$ is a geodesic if and only if
its coordinates $(x^i(t))$ satisfy
\be
\ddot{c}^i + 2 G^i (\dot{c}) =0.
\ee
$F$ is said to be {\it positively complete} (resp. {\it negatively complete}), if any geodesic on an open interval $(a, b)$ can be extended to a geodesic on $(a, \infty)$ (resp. $(-\infty, b)$). $F$ is said to be {\it complete} 
if it is positively and negatively complete.
There are  Finsler metrics which are positively complete, but not 
complete. See Example \ref{exBaoShen} below.

\bigskip
The notion of Riemann curvature for Riemann metrics can be extended to Finsler metrics/sprays. For a vector $y\in T_xM-\{0\}$,
the Riemann curvature ${\bf R}_y: T_xM\to T_xM$ is defined by
\[ {\bf R}_y (u) = R^i_k (y) u^k \; \pxi,\]
where
\be
 R^i_k (y):= 2 {\pa G^i\over \pa x^k}
- {\pa^2 G^i\over \pa x^j \pa y^k} y^j + 2 G^j {\pa^2 G^i\over \pa y^j \pa y^k} - {\pa G^i\over \pa y^j} {\pa G^j \over \pa y^k}.\label{Rik}
\ee

Take an arbitrary plane $P\subset T_xM$ (flag) and a non-zero vector $y\in P$ (flag pole),
the {\it flag curvature} $K(P, y)$ is defined by
\be
K(P, y) := { g_y({\bf R}_y (v), v) \over g_y(y,y) g_y(v,v)-g_y(v,y)g_y(v, y) }.
\ee
$F$ is said to be {\it scalar curvature} if for any nonzero vector $y\in T_xM$ and any
 flag $P\subset T_xM$, $x\in M$, with $y\in P$, 
$K(P, y) = \lambda(y)$ is independent of $P$, or
equivalently,
\[ {\bf R}_y = \lambda (y) F^2(y)\; \Big \{  I - g_y (y, \cdot ) y\Big \}, \ \ \ \ \ y\in T_xM, \ x\in M,\]
where $I:T_xM\to T_xM$ denotes the identity map and $g_y(y, \cdot) = -{1\over 2} [F^2]_{y^i} dx^i$. 
It is said to be {\it of constant curvature $\lambda$} if the above identity holds for the constant $\lambda$.

The trace of the Riemann curvature ${\bf R}_y$ is a scalar function ${\bf Ric}$ on 
$TM$
\be
{\bf Ric}(y):= {\rm trace \; of } \; {\bf R}_y .
\ee
${\bf Ric}$ is called the {\it Ricci curvature}.
Let $R = {1\over n-1}{\bf Ric}$.

A Finsler metric is called an {\it Einstein metric
with Einstein constant} $\lambda$ if 
\be
{\bf Ric}(y) = (n-1) \; \lambda \; F^2(y) .\label{Ricconstant1}
\ee
(\ref{Ricconstant1}) is simply denoted by ${\bf Ric}= (n-1) \lambda$ if no confusion is caused. 

\bigskip
We now consider pointwise projectively related Finsler metrics --- those having the same geodesics as set points. Given two Finsler metrics $F$ and $\tilde{F}$ on an $n$-dimensional manifold $M$, let ${\bf G}$ and $\tilde{\bf G}$ be the sprays induced by $F$ and $\tilde{F}$, respectively. It is easy to verify that 
\be
\tilde{G}^i = G^i + {\tilde{F}_{;k}y^k\over 2 \tilde{F}} \; y^i 
+ {\tilde{F}\over 2} \tilde{g}^{il} \Big \{ 
{\pa \tilde{F}_{;k}\over \pa y^l} y^k - \tilde{F}_{;l} \Big \}, \label{ghgh}
\ee
where $\tF_{;k}$ denotes the covariant derivatives of $\tF$ on $(M, F)$.
\be
 \tilde{F}_{;k} := {\pa \tilde{F}\over \pa x^k}
-{\pa G^l\over \pa y^k} {\pa \tilde{F}\over \pa y^l}.\label{covariant}
\ee
The identity (\ref{ghgh}) was first  established by A. Rapcs\'{a}k \cite{Rap}. By (\ref{ghgh}), Rapcs\'{a}k
 proved the following important lemma

\begin{lem} \label{Raplem}{\rm (Rapcs\'{a}k)} Let $(M, F)$ be a Finsler space. A Finsler metric $\tilde{F}$ is pointwise projective to $F$ if and only if 
\be
{\pa \tilde{F}_{;k}\over \pa y^l}  y^k - \tilde{F}_{;l} =0.\label{Rap1}
\ee
In this case, 
\be
 \tilde{G}^i = G^i + P y^i  \label{GiP}
\ee
with 
\be
P = { \tilde{F}_{;k} y^k \over 2 \tilde{F}} .\label{PF}
\ee 
\end{lem}

By Rapcs\'{a}k's lemma, we conclude that a Finsler metric $\tilde{F}$ on an open subset $\Omega\subset \R^n$ is a projective metric if and only if 
\be
{\pa \tilde{F}\over \pa x^k \pa y^l} y^k - {\pa \tilde{F}\over \pa x^l} = 0.\label{projRap}
\ee
One can verify that the Klein metric  $F_K$ and the Funk metrics $F_{\pm}$ in Section 1.1  satisfy (\ref{projRap}). Hence they are projective metrics. 

\bigskip

Let $F$ and $\tilde{F}$ be  Finsler metrics on an $n$-dimensional manifold $M$.
Assume that $\tilde{F}$ is pointwise projective to $F$, i.e., it satisfies (\ref{Rap1}).
Plugging (\ref{GiP}) into (\ref{Rik}) yields
\begin{eqnarray}
\tilde{\bf R}_y(u) & = & {\bf R}_y(u) + \Xi (y) \; u + \tau_y(u) \; y,\label{Rap6a}\\
\widetilde{\bf Ric}(y) & = & {\bf Ric} (y)
+ (n-1) \Xi(y),\label{Rap6}
\end{eqnarray}
where
\begin{eqnarray}
\Xi(y): & = &  P^2 - P_{;k}y^k , \label{Xi}\\
\tau_y (u) : & = & 3 \Big ( P_{;k} - {1\over 2} {\pa [P^2] \over \pa y^k}  \Big )u^k + {\pa \Xi\over \pa y^k}  u^k,\label{tau}
\end{eqnarray}
where $P_{;k}$ denote the covariant derivatives of $P$ on $(M, F)$ as defined in (\ref{covariant}) for $\tF$.
Using (\ref{PF}), one can express $\Xi(y)$ and $\tau_y$ in terms of $\tF$ and its covariant derivatives on $(M, F)$.
These formulas are given in   \cite{M1}\cite{MW}. 

We immediately obtain the following 
\begin{prop}\label{prop2.1} Let $(M, F)$ be a Finsler space of dimension $n$ and $\tilde{F}$ another Finsler metric on $M$. 

(i) Assume that ${\bf Ric} = (n-1) \lambda$. Then $\widetilde{\bf Ric} = (n-1)\tilde{\lambda}$ if and only if 
\be
 \Xi = \tilde{\lambda} \tilde{F}^2 - \lambda F^2. \label{Xishen}
\ee

(ii) Assume that ${\bf R}= \lambda$. Then $\tilde{\bf R} = \tilde{\lambda}$ if and only if 
(\ref{Xishen}) holds. 
\end{prop}

Proposition \ref{prop2.1}(ii) is proved in \cite{MW}.

\bigskip
There is a simple sufficient condition for $\tF$ being  of negative constant curvature and pointwise projective to $F$. 
\begin{prop}\label{prop2.2} Let $(M, F)$ be a Finsler space of dimension $n$ and $\tF$ another Finsler metric on $M$. Suppose that 
\be
 \tF_{;k} =  \mu {\pa [\tF^2]\over \pa y^k} , \label{FgFgg}
\ee 
where $\mu$ is a constant.
Then $\tF$ is pointwise projective to $F$ and 
\begin{eqnarray}
 \tilde{\bf R}_y (u)& = & {\bf R}_y (u) - {\mu^2}
\Big [  \tilde{F}^2(y) u - \tilde{g}_y (y, u ) y\Big ] ,\\
\widetilde{\bf Ric} (y) & = & {\bf Ric} (y) - (n-1) \mu^2 \tF^2(y),
\end{eqnarray}
where $\tilde{g}_y (y, u) := {1\over 2} {\pa [\tF^2]\over \pa y^k}  (y) u^k $.
Hence,

 (i) if $F$ is Ricci-flat($ {\bf Ric}=0$), then $\tF$ is an Einstein metric with $\widetilde{\bf Ric} = - (n-1)\mu^2 $;

(ii)  if $F$ is $R$-flat (${\bf R}=0$), then $\tF$ is of constant curvature with $\tilde{\bf R}= -\mu^2$.
\end{prop}
{\it Proof}: Differentiating (\ref{FgFgg}) with $y^l$, we obtain
\be
 {\pa \tF_{;k}\over \pa y^l}  = \mu {\pa^2 [\tF^2]\over \pa y^k \pa y^l}.\label{FgFggg}
\ee
Contracting (\ref{FgFggg}) with $y^k$ and using (\ref{FgFgg})
again, we obtain
\be
{\pa \tF_{;k} \over \pa y^l} y^k =\mu {\pa^2 [\tF^2]\over \pa y^k \pa y^l}y^k =   \mu {\pa [\tF^2]\over \pa y^l}= \tF_{;l} . 
\ee
By Lemma 2.1, we conclude that $\tF$ is pointwise projective to $F$.
Contracting (\ref{FgFgg}) with $y^k$ yields
\be
\tF_{;k}y^k = 2\mu \tF^2.\label{eq29}
\ee
Thus the function $P$ in (\ref{PF}) simplifies to
\be
 P  = { \tF_{;k} y^k \over 2 \tF} = \mu \tF. \label{P30}
\ee
Using (\ref{Xi}), (\ref{eq29}) and (\ref{P30}), we obtain
\begin{eqnarray*}
\Xi(y) & = & \mu^2 \tF^2 - \mu \tF_{; k} y^k\\
& = & \mu^2 \tF^2 - 2 \mu^2 \tF^2 = - \mu^2 \tF^2.
\end{eqnarray*}
Plugging (\ref{P30}) into (\ref{tau}) yields
\[ \tau_y(u) = \mu^2 \tilde{g}_y (y, u) .\]
This proves the proposition. 
\qed

\bigskip
The Funk metric $F_{\pm}$ on a strongly convex domain $\Omega\subset \R^n$  
satisfy
\be
 {\pa F_{\pm}\over \pa x^k} =  \pm {1\over 2} {\pa [F_{\pm}^2] \over \pa y^k}.\label{OkOk}
\ee
Thus $F_{\pm}$ are of constant curvature $-1/4$ by Proposition \ref{prop2.2}. Since $F_{-}$ and $F_{+}$ are projective metrics,
so is $F_H= {1\over 2}(F_{-}+F_{+})$. 
It follows from (\ref{OkOk}) that the projective factor $P$  of $F_H$ in (\ref{PF}) is given by
\[ P := { [F_H]_{;k}y^k\over 2 F_H} 
= {1\over 2} (F_{+} -F_{-}),\]
and 
\[ \Xi := P^2 - P_{;k}y^k = - [F_{H}]^2 .\]
Therefore the Hilbert metric
$F_H$ has constant curvature $-1$. 
This proof  is given  by Okada \cite{Ok}.

\section{Projectively Related Einstein Metrics}

Assume that  $F$ and $\tilde{F}$ are pointwise projectively related  Einstein metrics with
\[ {\bf Ric}(y) = (n-1)\lambda {F}^2(y) , \ \ \ \ \ \widetilde{\bf Ric}(y) = (n-1)\tilde{\lambda} \tilde{F}^2(y).\]
Then  (\ref{Rap1}) and (\ref{Xishen}) hold. It is much easier to work on (\ref{Xishen}) than (\ref{Rap1}).
Let us write (\ref{Xishen}) as follows.
\be
\tilde{\lambda} \tilde{F}^2 = \lambda F^2 
+ {3\over 4} \Big ( { \tilde{F}_{;k}y^k \over \tilde{F} } \Big )^2 
- { \tilde{F}_{;k;l}y^ky^l\over 2 \tilde{F}}.\label{Rap7}
\ee

Let $c(t)$ be an arbitrary unit speed geodesic in $(M, F)$
and 
\[ \tilde{F}(t):= \tilde{F}(\dot{c}(t)).\]
Observe that
\[ \tilde{F}'(t)=\tilde{F}_{;k}(\dot{c}(t)) \dot{x}^k(t), \ \ \ \ \tilde{F}''(t)= \tilde{F}_{;k;l}(\dot{c}(t)) \dot{x}^k(t)\dot{x}^l(t) .\]
Let 
\[ f(t):= {1\over \sqrt{\tilde{F}(t)}}.\]
(\ref{Rap7}) simplifies to
\be
f''(t) + \lambda f(t) = {\tilde{\lambda}\over f^3(t)}.\label{Rap9}
\ee
The equation (\ref{Rap9}) is solvable.

\bigskip
For simplicity, let 
\[
C:= {1\over 2} \Big (\lambda a^2 +\tilde{\lambda}/a^2 + b^2\Big ).\]
The  solution of (\ref{Rap9}) with 
\[ f(0)=a >0, \ \ \ \ f'(0)=b\not=0\]
 is  determined by
\be
\int^{f(t)}_{a} { s \over 
\sqrt{-\lambda s^4 + 2C  s^2 -\tilde{\lambda} } } ds = \pm t,\label{cc1}
\ee
where the sign $\pm$ in (\ref{cc1}) is same as that of $f'(0)=b$.
The solution with 
\[ f(0)=a >0, \ \ \ \ \ f'(0)=0\]
can be obtained by letting $b\to 0$.

Note that 
\be
-\lambda \Big ( a^2 - C/\lambda \Big )^2 + C^2/\lambda -\tilde{\lambda} 
= (ab)^2>0, \ \ \ \ {\rm if}\ \lambda \not=0  \label{eq23}
\ee
and
\be
-\lambda a^4 + 2C a^2 -\tilde{\lambda} = 
 (ab)^2 >0.\label{eq24}
\ee
Thus the integrand in (\ref{cc1}) is defined for $s$ close to $a$
and the maximal solution $f(t) >0$ exists on an interval $I$ containing $s=0$.

\section{$\lambda =1$}
 
In this section, we study the equation (\ref{cc1}) when $\lambda =1$. In this case
\[ C = {1\over 2} \Big ( a^2 +\tilde{\lambda}/a^2 +b^2 \Big ), 
\ \ \ \ C^2 -\tilde{\lambda} = \Big ( a^2 -C \Big )^2 + (ab)^2.\]

From (\ref{cc1}), we obtain
\be
f(t) = \sqrt{ (a^2-C) \cos (2t) + ab \sin(2t) + C}. \label{eqcase1}
\ee  
We use (\ref{eq23}) to rewrite (\ref{eqcase1}) in the following form
\be
f(t) = \sqrt{ \sqrt{C^2-\tilde{\lambda}} \; \sin \Big [ \sin^{-1} \Big ( {a^2-C\over \sqrt{C^2-\tilde{\lambda}}}\Big ) \pm 2t\Big ] +C},\label{eqcase1A}
\ee
where the sign $\pm$ in (\ref{eqcase1A}) is same as that of $f'(0)=b$ when $b\not=0$. Otherwise, the sign can be chosen arbitrarily.

\bigskip
\noindent 
{ \bf Case 1}: $\tilde{\lambda} =1$. In this case,
\[C= {1\over 2} (a^2+1/a^2+b^2) > 1, \ \ \ C^2-1 = (a^2-C)^2 + (ab)^2,\]
\[
 {|C|\over \sqrt{C^2-1}} > 1.\]
 Then  \[ f(t) = \sqrt{ \sqrt{C^2 -1} \sin \Big [\sin^{-1}
\Big ( {a^2-C\over \sqrt{C^2-1}} \Big ) \pm 2t   \Big ] +C}.\] 
Thus $f(t)$ is defined on
$I = (-\infty, \infty)$ and
for any $ r$, 
\[ \int_{r}^{r+\pi} {1\over f(t)^2} dt = \pi.\]

\bigskip
\noindent{\bf Case 2}: $\tilde{\lambda}= 0$. 
In this case, 
\[ C= {1\over 2} (a^2+b^2) >0, \ \ \ \ C^2 = (a^2-C)^2 + (ab)^2.\]
Then \[ f(t) =\sqrt{C} \sqrt{ \sin \Big [ \sin^{-1} \Big ( {a^2-C\over C} \Big ) \pm 2t \Big ] +1 }.\]
Thus $f(t)$ is defined on a bounded interval $I = (-\delta, \tau)$ and 
\[ \int_{-\delta}^0 {1\over f(t)^2} dt
= \infty = \int_0^{\tau} {1\over f(t)^2} dt .\]

\noindent{\bf Case 3}: $\tilde{\lambda} =-1$.
In this case, 
\[ C= {1\over 2} ( a^2 -1/a^2 +b^2), \ \ \ \ C^2+1 = (a^2-C)^2 + (ab)^2,\]
\[ {|C|\over \sqrt{C^2+1} } < 1.\]
Then
\[ f(t) = \sqrt{ \sqrt{C^2+1} \sin \Big [  \sin^{-1} \Big ( {a^2-C\over \sqrt{C^2+1} } \Big ) \pm 2t  \Big ] + C }.\]
Thus $f(t)$ is defined on a bounded interval  $I = (-\delta, \tau)$ and 
\[ \int_{-\delta}^0 {1\over f(t)^2} dt
= \infty = \int_0^{\tau} {1\over f(t)^2} dt .\]

\bigskip

From the above arguments, we obtain the following

\begin{prop} \label{lemRap0} Let $F$ and $\tilde{F}$ be Einstein metrics on an $n$-manifold $M$ with
\[ {\bf Ric} = n-1, \ \ \ \ \widetilde{\bf Ric} = (n-1)\tilde{\lambda}.\]
Then for any unit speed geodesic $c(t)$ of $F$, 
\be
 \tilde{F}(\dot{c}(t))
={2 \over 
\Big ( a^2-\tilde{\lambda}/a^2 - b^2 \Big ) 
\cos (2t) + 2 ab \sin (2t) + 
\Big ( a^2 +\tilde{\lambda}/a^2 + b^2 \Big )  }.\label{sol1}
\ee

(i) If $\tilde{\lambda} = 1$, then along any unit speed geodesic $c$ of $F$,
\[ \tilde{F}(\dot{c}(t)) 
= { 1\over \sqrt{C^2-1} \sin ( \theta  \pm 2t ) + C },\]
where $C>1$ and $\theta \in [-{\pi/2}, {\pi/2}]$.
Thus for any unit speed geodesic $c$ of $F$ with $F$-length
$L_F (c) =\pi$, the $\tilde{F}$-length $L_{\tilde{F}}(c) = \pi$. 

\smallskip

(ii) If $ \tilde{\lambda}=0$, then along any unit speed geodesic $c$ of $F$,
\[ \tilde{F} (\dot{c}(t))
= {1\over C  \sin ( \theta \pm 2t ) + C  },\]
where $C>0$ and  $\theta \in [-{\pi/2}, {\pi/2}]$.
Thus 
every  geodesic of $F$ has finite length. 

\smallskip

(iii) If $\tilde{\lambda}=-1$, the along any unit speed geodesic $c$ of $F$, 
\[ \tilde{F}(\dot{c}(t)) = {1\over \sqrt{C^2 +1} 
\sin ( \theta \pm 2t ) + C }, \]
where $C$ is a constant and  $\theta \in [-{\pi/2}, {\pi/2}]$.
Thus every  geodesic of $F$ has finite length. 
\end{prop}

\bigskip
Consider the following spherical metric on $\R^n$:
\be
F_S(y) :={ \sqrt{|y|^2 + \Big ( |x|^2 |y|^2 - \langle x, y \rangle^2 \Big )}\over  1 + |x|^2  },\label{exF_S} 
\ee
$F_S$ is of  constant curvature $1$ and  pointwise projective to the standard complete Euclidean metric $F_E$ on ${\rm R}^n$. More over,  all geodesics (straight lines) have $F_S$-length
$L_{F_S}(c) =\pi$. 

\section{$\lambda=0$}

In this section, we shall study the equation (\ref{cc1}) when $\lambda =0$. 
From (\ref{cc1}) we obtain
\be
f(t) = \sqrt{\Big ( a + bt\Big )^2 + \tilde{\lambda} \Big ({t\over a } \Big )^2 }. \label{eqcase2*}
\ee

\bigskip

\noindent{\bf Case 1}: $\tilde{\lambda} =1$. In this case, 
 Then 
\[ f(t) = \sqrt{ \Big ( a + bt \Big )^2 + \Big ( {t\over a }   \Big )^2} .\]
Thus $f(t)$ is defined on $I=(-\infty, \infty)$ and 
\[ \int_{-\infty}^{\infty} {1\over f(t)^2 } dt =\pi.\]

\bigskip

\noindent{\bf Case 2}: $\tilde{\lambda} =0$. In this case, 
\[ f(t) = a + bt .\]

(i)
If $b=0$, then 
\[ f(t) = a .\]
Thus $f(t)$ is defined on $I= (-\infty, \infty)$.

(ii)
If $b\not=0$, then 
\[ f(t) = a + bt.\]
In this case when $b >0$, $I = (-\delta, \infty)$ and 
\[ \int_{-\delta}^0 {1\over f(t)^2} dt =\infty \ \ \ 
{\rm and} \ \ \ \int_0^{\infty} {1\over f(t)^2} dt <\infty.\]
The case when $b <0$ is similar, so is omitted.

\bigskip

\noindent{\bf Case 3}: $\tilde{\lambda} = -1$. In this case,
\[ f(t) = \sqrt{ \Big ( a+bt \Big)^2 - \Big ( {t\over a }  \Big )^2  } .\]

(i) If  $ab = 1$, then 
\[ f(t) = \sqrt{ a^2 + 2 t }.\]
In this case, $ I = (-\delta, \infty)$ and 
\[ \int_{-\delta}^{0} {1\over f(t)^2 } dt = \infty = \int_0^{\infty} {1\over f(t)^2} dt.\]
The case  when $ab=-1$ is similar.

(ii) If $ -1 < ab < 1$, then $f(t)$ is defined on a bounded interval $(-\delta, \tau)$ and 
Clearly,
\[ \int_{-\delta}^0 {1\over f(t)^2} dt =\infty \ \ \ {\rm and} \ \ \
\int_0^{\tau} {1\over f(t)^2} dt =\infty.\]

(iii) If 
 $ ab > 1$, then $f(t)$ is defined on $( -\delta, \infty)$ and 
\[ \int_{-\delta}^{0} {1\over f(t)^2} dt =\infty \ \ \
{\rm and} \ \ \ \int_0^{\infty} {1\over f(t)^2} dt <\infty.\]
The case when $ab < -1$ is similar, so is omitted.

\bigskip

\begin{prop} \label{lemRap11}
Let $F$ and $\tilde{F}$ be Einstein metrics on an $n$-manifold $M$ with 
\[ {\bf Ric} = 0, \ \ \ \ \ \widetilde{\bf Ric} = (n-1)\tilde{\lambda}.\]
Assume that $F$ and $\tilde{F}$ are pointwise projectively related on $M$. 
Then for any unit speed geodesic $c(t)$ of $F$, 
\be
\tilde{F}(\dot{c}(t)) 
= {1\over \Big ( a+bt\Big )^2 + \tilde{\lambda}\Big ({t\over a} \Big )^2 } .
\ee

(i) If $\tilde{\lambda} =1$, then along any geodesic $c(t)$ of $F$,
\be
 \tilde{F} (\dot{c}(t)) =
{1\over \Big ( a+bt \Big )^2 + \Big ( {t\over a} \Big )^2 }.\label{32a}
\ee
Thus 
for any geodesic $c$ of $F$, the $\tilde{F}$-length  $ L_{\tilde{F}} (c) \leq  \pi$. Equality holds when 
$F$ is complete.

(ii) If $\tilde{\lambda}=0$, then along any unit speed geodesic $c(t)$ of $F$,
\[ \tilde{F}(\dot{c}(t)) = {1\over (a+bt)^2 }.\]

(iia) If a unit speed geodesic $c$ of $F$ is defined on $(-\infty, \infty)$, then 
$ \tilde{F}(\dot{c}(t)) = 1/a^2.$

(iib) If a unit speed geodesic $c$ of $F$ is defined on $[0, \infty)$, then 
it has finite $\tilde{F}$-length unless $\tilde{F}(\dot{c}(t)) = 1/a^2.$

(iib) If a unit speed geodesic $c$ of $F$ is defined on $(-\infty, 0]$, then 
it has finite $\tilde{F}$-length unless $\tilde{F}(\dot{c}(t)) = 1/a^2.$

\noindent
Therefore,  $F$ is complete if and only if  $\tilde{F}$ is complete. In this case,  along any geodesic $c$
\[ {{F}(\dot{c}(t))\over \tilde{F}(\dot{c}(t))} = constant.\]

(iii) If $\tilde{\lambda}= -1$, then along any unit speed geodesic $c$ of $F$,
\be
 \tilde{F}(\dot{c}(t)) 
= {1\over \Big ( a + bt \Big )^2 - \Big ( {t \over a } \Big )^2 }.\label{32b}
\ee
Thus no geodesic of $F$ is defined on $(-\infty, \infty)$.

(iiia) If a unit speed geodesic $c$ of $F$ is defined on $[0, \infty)$, then
it has finite $\tilde{F}$-length unless 
\be
 \tilde{F}(\dot{c}(t)) = {1\over a^2 + 2t }.\label{32c}
\ee

(iiib) If a unit speed geodesic $c$ of $F$ is defined on $(-\infty, 0]$, then 
it has finite $\tilde{F}$-length unless
\be
 \tilde{F}(\dot{c}(t)) = {1\over a^2 - 2t }.\label{32d}
\ee

\end{prop}

\bigskip

For the spherical metric  $F_S$ in (\ref{exF_S}),
the geodesics of $F_S$ are straight lines in $\R^n$. Thus it is pointwise projective to the standard Euclidean metric $F_E$ on ${\rm R}^n$. It is easy to verify that all lines in ${\rm R}^n$  have length of $\pi$ with respect to $F_S$. 
This is also guaranteed by Proposition \ref{lemRap11} (i).

\bigskip
\section{$\lambda = -1$}
In this section, we shall study the equation (\ref{cc1}) when
${\lambda}=-1$. In this case,
\[ C = {1\over 2} \Big ( -a^2 +\tilde{\lambda}/a^2 + b^2 \Big ), 
\ \ \ \ \Big ( a^2 + C \Big )^2 = C^2 +\tilde{\lambda} +(ab)^2.\]
From (\ref{cc1}), we obtain
\be
f(t)= \sqrt{(a^2+C) \cosh (2t) + ab \sinh (2t) - C}.\label{eqcasexA}
\ee
We use (\ref{eq23}) to rewrite (\ref{eqcasexA}) as follows
\be
f(t)
=\cases{\sqrt{ \sqrt{C^2 +\tilde{\lambda}} \; \cosh \Big [\cosh^{-1} \Big ( { a^2 + C\over \sqrt{ C^2 +\tilde{\lambda}} } \Big ) \pm 2 t
\Big  ]- C } & if $C^2 +\tilde{\lambda} >0$\cr\\
 \sqrt{ e^{\pm 2t} \Big ( a^2 + C \Big ) -C}& if $C^2 +\tilde{\lambda} =0$\cr\\
 \sqrt{ \sqrt{-C^2 -\tilde{\lambda}} \; \sinh \Big [\sinh^{-1} \Big ( { a^2 + C\over \sqrt{ -C^2 -\tilde{\lambda}} } \Big ) \pm 2 t
\Big  ]- C }& if $C^2 +\tilde{\lambda } < 0$} \label{eqcase3*}
\ee
The sign $\pm$ in (\ref{eqcase3*}) is same as that  of $f'(0)=b\not=0$.

 We divide this case into several cases.

\bigskip

\noindent{\bf Case 1}: $\tilde{\lambda} =1$.
In this case, 
\[ C = {1\over 2} \Big ( -a^2 + 1/a^2 +b^2 \Big ), \ \ \ \  (a^2+C)^2 = C^2 +1 + (ab)^2,\]  
\[{ |C|\over \sqrt{C^2+1}}  < 1.\]
Then 
\[ f(t) = \sqrt{ \sqrt{C^2 +1} \; \cosh \Big [\cosh^{-1} \Big ( { a^2 + C\over \sqrt{ C^2 +1} } \Big ) \pm 2 t
\Big  ]- C }.\]
Thus $f(t)$ is defined on $I=(-\infty,\infty)$ and
\[ \int_{-\infty}^{\infty} {1\over f(t)^2} dt <\infty.\]

\bigskip
\noindent{\bf Case 2}: $\tilde{\lambda} =0$. In this case,
\[  C= {1\over 2} ( -a^2 + b^2 ), \ \ \ \ (a^2+C)^2 = C^2 + (ab)^2.\]
Then 
\[ f(t) = \cases { \sqrt{2|C|} \cosh \Big [ {1\over 2} \cosh^{-1} \Big ( {a^2\over |C|} -1 \Big )   \pm t\Big ]  & if $ C<0$,\cr\\
\sqrt{2C} \sinh \Big [{1\over 2} \cosh^{-1} \Big ( {a^2\over C}+1 \Big )    \pm t \Big ]     & if $C >0$,\cr\\
a e^{\pm t} & if $C=0$  .\cr }
\]

(i) If $C<0$, then $f(t)$ is defined on $I= (-\infty, \infty)$ and
\[ \int_{-\infty}^{\infty} {1\over f(t)^2} dt <\infty.\]

(ii) If $C>0$, then $f(t)$ is defined on  either $I=(-\delta, \infty)$  or $I = (-\infty, \tau)$ . Assume that $I= (-\delta, \infty)$. Then 
\[ \int_{-\delta}^0 {1\over f(t)^2} dt =\infty \ \ \ 
{\rm and}\ \ \ \int_0^{\infty} {1\over f(t)^2} dt <\infty.\]
The case when $I=(-\infty, \tau)$ is similar, so is omitted.

(iii) If $C=0$, then $b\not=0$ and  $f(t)$ is defined on $I = (-\infty, \infty)$. Assume that  $ b>0$. Then  
\[ \int_{-\infty}^0 {1\over f(t)^2} dt =\infty \ \ \ 
{\rm and} \ \ \ \int_0^{\infty} {1\over f(t)^2} dt < \infty.\]
The case when $b <0$ is similar, so is omitted.

\bigskip
\noindent{\bf Case 3}: 
$\tilde{\lambda} = -1$. 
In this case, 
\[ C= {1\over 2} \Big ( -a^2 - 1/a^2 + b^2 \Big ), \ \ \ \ (a^2+C)^2 = C^2-1 + (ab)^2.\]

{(i) } $ C^2 > 1$. In this case,
\[ {|C|\over \sqrt{ C^2-1}} >1 .\]
Then 
\[ f(t) = \sqrt{ \sqrt{C^2-1} \cosh \Big [ \cosh^{-1} \Big ({a^2+C\over \sqrt{C^2-1}}    \Big )  \pm 2t\Big ] - C}.\]

(ia)
If $C >1$, then $f(t)$ is defined on $I=(-\delta, \infty)$ and
\[ \int_{-\delta}^0 {1\over f(t)^2 } dt =\infty \ \ \
{\rm and} \ \ \ \int_0^{\infty} {1\over f(t)^2} dt <\infty.\]

(ib)If $C < -1$, then $f(t)$ is defined on $I= (-\infty, \infty)$ and
\[ \int_{-\infty}^{\infty} {1\over f(t)^2} dt <\infty.\]

{(ii)} $ C^2 <1 $. Then 
\[ f(t) = \sqrt{ \sqrt{1-C^2} \sinh \Big [ \sinh^{-1} \Big ( { a^2 +C \over 1-C^2 }  \Big )   \pm 2t \Big ] - C }. \]
In this case, $f(t)$ is defined on either $I = (-\delta, \infty)$ or $I=(-\infty, \tau)$. Assume that $ I= (-\delta, \infty)$. Then 
\[ \int_{-\delta}^0 {1\over f(t)^2} dt =\infty
\ \ \ {\rm and} \ \ \ \int_0^{\infty} {1\over f(t)^2} dt <\infty.\]
The case when $I= (-\infty, \tau)$ is similar, so is omitted.

{(iii) } $C^2 =1$. 

(iiia)
If $C= 1$, then 
\[ f(t) = \sqrt{ e^{\pm 2t} ( a^2+1) - 1 }.\]
Thus $f(t)$ is defined on either  $I= (-\delta, \infty)$ 
or $I= (-\infty, \tau)$. Assume that $ I = (-\delta, \infty)$. Then
\[ \int_{-\delta}^0 {1\over f(t)^2} dt =\infty \ \ \
{\rm and} \ \ \ \int_0^{\infty} {1\over f(t)^2} dt <\infty.\]
The case when $I= (-\infty, \tau)$ is similar, so is omitted.

(iiib)
If $C= -1$, then 
\[ f(t) = \sqrt{ e^{\pm 2t} ( a^2-1) + 1 }.\]
If $ a >1$, then $b\not=0$ and  $f(t)$ is defined on $I=(-\infty, \infty)$. Assume that $ b > 0$, then 
\[ \int_{-\infty}^0 {1\over f(t)^2 } dt = \infty 
\ \ \  {\rm and} \ \ \ \int_0^{\infty} {1\over f(t)^2} dt <\infty.\]
The case when $ b<0$ is similar, so is omitted.

If $ 0 < a < 1$, then  $f(t)$ is defined on either  $I= (-\delta, \infty)$ or  $I= (-\infty, \tau)$. Assume that $I= (-\delta, \infty)$. Then
\[  \int_{-\delta}^0 {1\over f(t)^2} dt =\infty
\ \ \ {\rm and} \ \ \ \int_0^{\infty} {1\over f(t)^2} dt  < \infty.\]
The case when $ b > 0$ is similar, so is omitted. 

If $a =1$, then $b=0$. Then 
\[ f(t) = 1.\]
In this case,  $f(t)$ is defined on $I= (-\infty, \infty)$.

\bigskip

\begin{prop}\label{lemRap1} Let $F$ and $\tilde{F}$ be Einstein metrics on an $n$-manifold 
$M$ with 
\[ {\bf Ric} = -(n-1), \ \ \ \ \ \widetilde{\bf Ric} = (n-1)\tilde{\lambda} .\]
Assume that $F$ and $\tilde{F}$ are pointwise projectively related. Then for any geodesic of $c(t)$ of $F$, 
\be
\tilde{F}(\dot{c}(t)) = {1\over 
\Big ( a^2 +\tilde{\lambda}/a^2 +b^2 \Big  )
\cosh (2t) + 2 ab \sinh (2t) -  \Big ( -a^2 +\tilde{\lambda}/a^2 +b^2 \Big )}.
\ee

(i)  $\tilde{\lambda} = 1$.  In this case,   any geodesic of $\tilde{F}$ has finite length. Hence $\tilde{F}$ is neither positively complete, nor negatively complete. 

\smallskip

(ii)  $\tilde{\lambda}=0$. In this case, no geodesic of $\tilde{F}$ is defined on $(-\infty, \infty)$.

(iia) If a unit speed geodesic $c$ of $F$ is defined on $[0, \infty)$, then it has finite $\tilde{F}$-length unless 
\be
 \tilde{F}(\dot{c}) = \Big (  {e^t \over a} \Big )^2.\label{+FF}
\ee

(iib) If a unit speed geodesic $c$ of $F$ is defined on $(-\infty, 0]$, then it has finite $\tilde{F}$-length unless
\be
 \tilde{F} (\dot{c}(t)) = \Big ( {e^{-t}\over a } \Big )^2 .\label{-FF}
\ee

\smallskip

(iii) $\tilde{\lambda}=-1$. In this case, if both  $F$ and $\tilde{F}$ are complete, then 
\[ F =\tilde{F}.\]

(iiia) If a unit speed geodesic $c$ of $F$ is defined on $[0, \infty)$, then 
it has finite $\tilde{F}$-length unless 
\be
 \tilde{F}(\dot{c}(t)) = { 1\over e^{-2t} (a^2 -1) + 1 }, \ \ \ \ \ ( a \geq 1).\label{+F}
\ee

(iiib) If a unit speed geodesic $c$ of $F$ is defined on $(-\infty, 0]$, then 
it has finite $\tilde{F}$-length unless
\be
 \tilde{F}(\dot{c}(t)) = { 1\over e^{2t} (a^2 -1) + 1 }, \ \ \ \ \ ( a \geq 1).\label{-F}
\ee

\end{prop}

\section{Examples}

Below are some interesting examples. All the metrics are projective Finsler metrics of constant curvature on a strongly convex domain in 
the Euclidean space. 

\begin{ex}\label{exCase2a} {\rm 
The standard metric on the upper/lower semi-sphere ${\rm S}^n_{\pm}$ can 
be pulled back to the spherical metric $F_S$ on ${\rm R}^n$ by a diffeomorphism $\varphi_{\pm}: \R^n \to {\rm S}^n_{\pm}$,
\be
 \varphi_{\pm}(x):= \Big ( {x \over \sqrt{1+|x|^2 }}, \; {\pm 1\over \sqrt{1+|x|^2 } }\Big ).\label{phipm}
\ee
The formula of $F_S$ is given in (\ref{exF_S}). 
$F_S$ has positive constant curvature $=1$ and it  is pointwise projective to the standard flat metric $F_E(y)=|y|$ on $\R^n$.
Take an arbitrary geodesic $c(t) = x + ty$ in $(\R^n, F)$. Then
\begin{eqnarray}
F_S(\dot{c}(t)) & = & 
{ \sqrt{|y|^2 + \Big ( |x|^2 |y|^2 - \langle x, y \rangle^2 \Big )}\over  1 + |x|^2+ 2 \langle x, y \rangle t + |y|^2 t^2   }\\
&=&  { 1\over \Big ( a + bt\Big )^2 + \Big ( {t\over a}  \Big )^2 },
\end{eqnarray}
where 
\begin{eqnarray*}
 a & = &  { \sqrt{1+|x|^2} \over \Big [ |y|^2 + \Big ( |x|^2 |y|^2 - \langle x, y \rangle^2\Big ]^{1/4}}\\
 b & = &  {\langle x, y \rangle 
\over \sqrt{ 1+|x|^2} \Big [ |y|^2 + \Big ( |x|^2 |y|^2 - \langle x, y \rangle^2\Big ]^{1/4} }\end{eqnarray*}
 Thus $F_S(\dot{c}(t))$ is in the form (\ref{32a}).}
\end{ex}

\begin{ex}\label{Bryantmetric}{\rm 
Deforming  the spherical metric $F_S$ yields some interesting 
Finsler metrics. Let 
\begin{eqnarray*}
A_{\e}(y) : & = &|x|^2 |y|^2 - \langle x, y\rangle^2 + \e \; |y|^2 + { 2 (1-\e^2)\langle x, y\rangle^2 \over |x|^4 + 2 \e |x|^2 +1 } ,\\
B_{\e}(y): & = & 
\Big ( |x|^2 |y|^2 - \langle x, y\rangle^2\Big )^2 + 2 \e\; 
\Big ( |x|^2 |y|^2 - \langle x, y\rangle^2\Big ) |y|^2 + |y|^4 .
\end{eqnarray*}
For $0 < \e \leq 1$, define
\be
 F_{\e}(y) := \sqrt{A_{\e}(y) +\sqrt{B_{\e}(y)}\over 2\Big (  |x|^4 +2\e |x|^2 +1 \Big )} 
+  { \sqrt{1-\e^2} \langle x, y \rangle \over 
 |x|^4 + 2\e \; |x|^2 + 1 }, \ \ \ y\in T_x\R^n.\label{FFBB}
\ee
$F_{\e}$ is a family of Finsler metrics on $\R^n$. 
Note that
$F_1=F_S$ is just the spherical metric in (\ref{exF_S}).
Using (\ref{phipm}), one can pull $F_{\e}$ onto 
${\rm S}^n$. The pull-back metrics  on ${\rm S}^n$ are the natural generalization of Bryant metrics on ${\rm S}^2$  \cite{Br1}\cite{Br2}.

Assume that $F_{\e}$ are of constant curvature $1$ and pointwise projective to the Euclidean metric $F_E$ on $\R^n$. Then 
 for any $c(t)=x+ty$, there are constants $a>0$ and $-\infty < b < \infty$ such that
\[ F_{\e}(\dot{c}(t)) = {1\over \Big (a +bt\Big)^2 + \Big ( {t\over a}  \Big )^2 }.\]
The  constants $a$ and $b$ must be given by
\[ a= {1\over \sqrt{F_{\e}(y)}}, \ \ \ \
b = - { y^i \over 
\sqrt{F_{\e}(y)}} {\pa \over \pa x^i} \Big [ \ln \sqrt{ F_{\e} (y) } \Big ].\]
The proof maybe need a faster computer.
}
\end{ex}

\bigskip
\begin{ex}\label{ex4.2}{\rm 
 The Klein metric $F_K$ in (\ref{ccHilbert}) is  Riemannian.
It is complete with  constant curvature $-1$ and  pointwise projective to the standard (incomplete) Euclidean metric $F_E(y)=|y|$ on ${\rm B}^n$.
Take an arbitrary  geodesic $c(t) = x + ty$ in $({\rm B}^n, F)$.
Then
\begin{eqnarray}
 F_K(\dot{c}(t)) & = & {\sqrt{ |y|^2 - \Big ( |x|^2 |y|^2 - \langle x, y \rangle^2 \Big ) }
\over 1- |x|^2- 2 \langle x, y \rangle t - |y|^2 t^2 }\nonumber\\
& = & {1\over \Big ( a+bt \Big )^2 - \Big (  {t\over a } \Big )^2 },\label{dcdcdc}
\end{eqnarray}
where 
\begin{eqnarray*}
 a & = &  {\sqrt{1-|x|^2} \over \Big [ |y|^2 - 
\Big ( |x|^2 |y|^2 - \langle x, y \rangle^2\Big ) \Big ]^{1/4} } \\ 
b & = & - {\langle x, y \rangle \over\sqrt{1-|x|^2}  \Big [ |y|^2 - \Big ( |x|^2 |y|^2 - \langle x, y \rangle^2\Big ) \Big ]^{1/4} }.
\end{eqnarray*}
Thus $F_K(\dot{c}(t))$ are in the form (\ref{32b}).
}
\end{ex}

\bigskip

\begin{ex}\label{exBaoShen}{\rm (Funk metrics)
Let 
$\Omega$ be a strongly  convex  bounded domain in ${\rm R}^n$. 
For  $0\not=y\in T_x\Omega \approx {\rm R}^n$,
let $ F_{-}(y) >0 $ and $F_{+}(y) >0$ be given by
\be
z_{-} = x - {y\over F_{-}(y)}, \ \ \ \ \ 
z_{+} = x + {y\over F_{+}(y)},\label{Funkpm}
\ee
where $z_{-}, z_{+}$ are the intersection points of the line $\ell(t):=x+ty$ with $\pa \Omega$ such that $z_{+}-z_{-}$ is in the same direction as $y$. 
$F_{\pm}$ are called the pair of {\it Funk metrics} on $\Omega$. Note that
$ F_{+}(-y) = F_{-}(y).$ More over, they satisfy
\[ {\pa F_{\pm}\over \pa x^k} = \pm {1\over 2}{\pa [F_{\pm} ]\over \pa y^k}.\]
According to Proposition \ref{prop2.1}, $F_{\pm}$ are of constant curvature $-1/4$ and pointwise projective to the 
standard Euclidean metric $F_E$ on $\Omega$. This simple proof is due to T. Okada \cite{Ok}.

Fix $y\in T_x\Omega$ and $t $ such that $c(t)=x+ty \in \Omega$. 
From the definition of $F_{-}$ and $F_{+}$, we have
\[ z_{-} = x - {y\over F_{-}(y)} = x+ ty - {y\over F_{-}( \dot{c}(t)) }.\]
\[ z_{+} = x + {y\over F_{+}(y)} = x+ ty + {y\over F_{+}( \dot{c}(t)) }.\]
Then we obtain
\begin{eqnarray}
F_{-} (\dot{c}(t)) & = & {F_{-}(y) \over 1+ F_{-}(y) t },\label{F1}\\
F_{+}(\dot{c}(t)) & = & {F_{+}(y)\over 1 - F_{+}(y) t }.\label{F2}
\end{eqnarray}
Thus
\[ {1\over 2} F_{\pm}(\dot{c}(t))
= {1\over a^2 \mp 2 t},\]
where $ a^2 = 2/F_{\pm}(y)$. Thus ${1\over 2} F_{\pm}(\dot{c}(t))$ are in the form (\ref{32b}) with $ab=\mp 1$. This is also guaranteed by Proposition \ref{lemRap11}, because that  ${1\over 2}F_{\pm}$ has constant curvature $-1$.

}
\end{ex}

\bigskip

\begin{ex}\label{exHM}
{\rm 
Let  $\Omega$  be a strongly convex bounded domain in ${\rm R}^n$.
Let $F_{\pm}$ denote the Funk metrics on $\Omega$ defined in (\ref{Funkpm}). 
 Define
\be
F_H(y) := {1\over 2}  \Big ( F_{-}(y) + F_{+}(y) \Big ).\label{FFFFHHHH}
\ee
$F_H$ is  called the {\it Hilbert metric} on $\Omega$.
The Hilbert is  of constant curvature $-1$ and pointwise projective to 
the Euclidean metric $F_E$ on $\Omega$. See \cite{Bu} \cite{BuKe} \cite{Fk} \cite{Ok}.

It follows from  Proposition \ref{lemRap11} that   along any geodesic $c(t)=x+ty$ of $F_E$,
$F_H(\dot{c}(t))$ should satisfy (\ref{32b}). Let us verify this necessary condition directly. 
From (\ref{F1}) and (\ref{F2}), we obtain 
\begin{eqnarray*}
F_H(\dot{c}(t)) & = & {F_{-}(y)+ F_{+}(y) \over 2 \Big ( 1+ F_{-}(y) t  \Big ) \Big (  1 - F_{+}(y) t  \Big ) }\\
& = & {1\over \Big ( a+bt \Big )^2 - \Big ( {t\over a } \Big )^2 },
\end{eqnarray*}
where 
\[ a = {\sqrt{2}\over \sqrt{F_{-}(y)+F_{+}(y)}}, \ \ \ \  b= {F_{-}(y)-F_{+}(y)
\over \sqrt{2}\sqrt{F_{-}(y)+F_{+}(y)}}.\]
Yes, $F_{H}(\dot{c}(t))$ satisfies (\ref{32b}).

\bigskip

Let $F=F_H$ and $\tilde{F}=F_{\pm}$. Then along any unit speed geodesic $c(t)$ of $F$,
$\tilde{F}_{+}(\dot{c}(t))$ satisfies (\ref{+F}) and $\tilde{F}_{-}(\dot{c}(t))$ satisfies (\ref{-F}).

 }
\end{ex}

\bigskip

\begin{ex}{\rm 
Let $\Omega$ denote the domain above the graph $ x^n = \sum_{a=1}^{n-1} (x^a)^2 $ in $\R^n$. Define $\tilde{F}: T\Omega \to [0, \infty)$ by
\be
\tilde{F} (y) := { \sqrt{\Big ( y^n - 2 \sum_{a=1}^{n-1} x^a y^a \Big )^2 
+ 4 \Big ( x^n - \sum_{a=1}^{n-1} (x^a)^2 \Big ) \sum_{a=1}^{n-1}
(y^a)^2} \over 2 \Big  ( x^n - \sum_{a=1}^n (x^a)^2 \Big ) }.
\ee 
This Riemann metric has constant curvature $-1$.
Take an arbitrary geodesic $c(t) = x + ty$ in $(\Omega, F)$, $\tilde{F}(\dot{c}(t))$ 
must be in the form (\ref{32b}).
Note that every geodesic $c(t)=x+ty$ must intersect $\pa \Omega$ on both sides unless
$ y = (0, \cdots, 0, y^n)$
with $y^n >0$. When $y = (0, \cdots, 0, y^n)$ with $y^n >0$,
the geodesic 
$c(t)=x+ty$ of $F$ is defined on 
$(-\delta, \infty)$ and the $\tilde{F}$-length of $c$ over  $[0, \infty)$ is finite. Moreover, 
\[  \tilde{F}(\dot{c}(t)) = { 1\over a^2 + 2 t },\]
where 
\[ a^2 = 2 { x^n- \sum_{a=1}^{n-1} (x^a)^2 \over y^n}.\] 
The last statement  is also implied by  Proposition \ref{lemRap11} (iiia).

}
\end{ex}

\bigskip
Finally we ask the following

\noindent
{\bf Open Problem}: Are there non-trivial  positively/negatively complete Ricci-flat metrics on an open subset $\Omega \subset \R^n$ ? If any, they must be $R$-flat (i.e., ${\bf R}=0$).

\end{document}